# ON CALCULATION OF CANONICAL DECOMPOSITION OF TENSOR VIA THE GRID OF LOCAL DISCREPANCIES

**A.K. Alekseev** [*], **A.E. Bondarev** [**], **Yu.S. Pyatakova** [*]

([*] *141070 Korolev, Lenin str. 4a, Korolev Rocket and Space Corporation Energia,*
[**] *125047 Moscow, Miusskaya str. 4, Keldysh Institute of Applied Mathematics RAS)*

*e-mail:* aleksey.k.alekseev@gmail.com, bond@keldysh.ru

**Abstract.** The method for calculation of the canonical decomposition that approximates a tensor of high order is considered, which requires moderate computational resources. It is based on the replacement of the approximation error norm (global discrepancy functional) by the grid of local functionals (discrepancies computed on hyperplanes). The point of the global functional minimum in the space of the canonical decomposition cores is determined by the set of the stationary points of local functionals. In result, the estimation of the cores of the canonical decomposition is possible using Newton method applied point-wisely along coordinates and nodes. The discrepancies on the hyperplanes are calculated using Monte-Carlo method.

Numerical tests on the approximation of sixth order tensors confirms the efficiency of the proposed approach.

## INTRODUCTION

At the operations with high order tensors, the application of the tensor decompositions (canonical decomposition [1,2] or tensor train [3]) enables the radical decrease of requirements both for the computer memory and the time of computation. Corresponding technique is developed in [3] for the tensor train and in [4] for the canonical decomposition. The paper [5] considers the approximation of a canonical decomposition by another one of the lesser rank via the Newton method that improves the set of algorithms provided by [4] for a reduction of the rank of decomposition.

However, the initial stage (the calculation of the tensor decompositions' coefficients) for the given tensor (approximation of a dense high order tensor) remains the complex problem from the computational viewpoint.

As a rule, the method of the alternating least squares (ALS) [9,10,11] is applied for the calculation of the coefficients of canonical decomposition. Conventionally, it is realized via the Khatri-Rao product [7] that uses the matricized tensor. A calculation of the Khatri-Rao product is computationally expensive and, for a tensor with an order higher 3, very difficult in realization. So, the attempts to construct alternative methods, such as provided in [11], are quite natural. Besides ALS the singular decomposition [3], also computationally expensive, is used for the calculation of the tensor train coefficients.

Most of the above mentioned methods in some way apply the complete tensor (usually in matricized form) or (and) its convolution on hyperplane that causes the exponential growth of required computational resources at increase the order of the tensor. This circumstance causes problems from the standpoints both of the memory and the computational time for the tensors of the high order.

The cross-approximation [6] is an exception that provides for the tensor train a linear on the tensor order growth of the time of computation and the computer memory.

In the present paper we consider for the canonical decomposition the approach that also provides linear growth of the required computational resources at increase of the tensor order and does not apply the total tensor and its' convolutions. This approach is based on the application of the Monte-Carlo method at approximation of multidimensional integrals (convolutions) at calculation of discrepancies and the replacement of the global functional (discrepancy of the tensor and its approximation) by the grid of local functionals (discrepancies on hyperplanes). This replacement is feasible since the gradient of the global functional may be computed by the derivatives of the local functionals.

The realizations of this approach using the Newton method and the steepest descent method are compared. The numerical tests confirmed the efficiency of both methods for the six order tensor approximation with the significant advantage of the Newton method regarding the computational time and accuracy.

## 1. METHODS FOR CANONICAL DECOMPOSITION CALCULATION

Let's consider the numerical algorithms for calculation of cores (factor-matrices) $Q_{\alpha,i}^j$ of canonical decomposition that approximates a high order tensor

$$A_{i_1,i_2...i_d} = \sum_{\alpha=1}^{r} Q_{\alpha,i_1}^1 Q_{\alpha,i_2}^2 ... Q_{\alpha,i_d}^d \tag{1}$$

(in our case obtained using values of a multidimensional function in nodes of regular grid). Herein indices $j = 1...d$ denotes directions (coordinates), indices $i_j$ mark the nodes along direction $j$, $\alpha$ is the number of the layer of core, $d$ is the dimension of the space (order of the tensor), $N_j$ is the number of nodes along single direction, as a rule, considered universal $N_j = N$, $r$ is the rank of canonical decomposition.

The one-dimensional on coordinates structure of the canonical decomposition cores in a fact dictates the application of the alternating direction methods. In our case, the alternating least squares (ALS) method [10,11] enables to minimize the functional of the discrepancy of exact data and approximation

$$\varepsilon(Q^1, Q^2, ..., Q^d) = \left\| A - \sum_{\alpha=1}^{r} Q_\alpha^1 Q_\alpha^2 ... Q_\alpha^d \right\|_{L_2}^2 / 2 \tag{2}$$

over one core at fixed others and to obtain the values of cores in a form

$$Q^k = \arg\min_{Q^k} \left\| A - \sum_{\alpha=1}^{r} Q_\alpha^1 Q_\alpha^2 ... Q_\alpha^d \right\|^2 / 2, k = 1...d \ . \tag{3}$$

This statement enables different ways of realization including solving the Euler-Lagrange equations [9], gradient optimization [12], Newton method [11,13,14,15], coordinate descent [16,17]. In most cases the Khatri-Rao product (marked as .) is used [7,9,12]. The cores (factor-matrices) $Q_k$ are determined by subsequent solution of the following problem

$$Q^k = \arg\min_{Q^k} \left\| A_{(k)} - Q^k (Q^1 . ... . Q^{k-1} . Q^{k+1} . ... . Q^d)^T \right\|^2, \ k = 1...d \ . \tag{4}$$

Here $A_{(k)}$ is the mode-k matricization of the tensor, $Q^k(Q^1 \cdot \ldots \cdot Q^{k-1} \cdot Q^{k+1} \cdot \ldots \cdot Q^d)^T$ is an auxiliary matrix.

In order to determine the minimum of (4) the following expression is used

$$Q^k = A_{(k)}((Q^1 \cdot \ldots \cdot Q^{k-1} \cdot Q^{k+1} \cdot \ldots \cdot Q^d)^T)^+, \tag{5}$$

where "+" notes a pseudoinverse matrix. This expression is commonly noted as matricized tensor times Khatri-Rao product, (MTTKRP) [8]. The matricization $A_{(k)}$ of the tensor used in expression (5) requires the same memory as the tensor. The inversion of matrix $(Q^1 \cdot \ldots \cdot Q^{k-1} \cdot Q^{k+1} \cdot \ldots \cdot Q^d)^T$ is also necessary. This matrix requires a bit less memory if compare with tensor. By these reasons the operation (5) is highly expensive both from the viewpoint of required memory and from the viewpoint of the time of computations.

So, the search for alternative algorithms is expedient. The most efficient (as far as we know) method is presented in [11], where the variant of ALS based on the local inversion of small ($r^2$) Hessians and the variant of the Newton method for a minimization of linearized discrepancy (inversion of Hessians $(d \times r^2)^2$) are used. These algorithms significantly reduce the required memory. Nevertheless, the summation over all indices of the tensor except one (convolution on the hyperplane) is used that causes high computational time for the tensors of high order.

In general, all known to authors methods both for a calculation of the canonical decomposition and for a calculation of the tensor train use the marticized tensor (in best case its rows and columns) and (or) convolutions of the tensor that causes high requirements to the memory and high time of computation for tensors of the order above three.

Herein we consider the method, which is significantly more inexpensive from the computational viewpoint, for the estimation of canonical decomposition cores. The method is based on gradients and Hessians for a grid of discrepancies computed on hyperplanes using Monte-Carlo method. As the test cases we consider the approximation by the canonical decomposition of the tensor defined by the values of multidimensional function in nodes of the grid.

## 2. ESTIMATION OF CANONICAL DECOMPOSITION IN VARIATIONAL STATEMENT FOR GLOBAL DISCREPANCY

The quality of the approximation of a tensor by the canonical decomposition may be estimated via the norm of the discrepancy $\widetilde{\varepsilon}_\Sigma = \left\| f_{i_1 \ldots i_d} - \widetilde{f}_{i_1 \ldots i_d} \right\|_{L_2}^2$, which (formally) may be computed by summation over all nodes

$$\widetilde{\varepsilon}_\Sigma = \sum_{i_s, s=1 \ldots d} \widetilde{R}_{i_1 \ldots i_d} \cdot \widetilde{R}_{i_1 \ldots i_d} / (2 \cdot N^d), \tag{6}$$

where

$$\widetilde{R}_{i_1 \ldots i_d} = \sum_{\alpha=1}^{r} \prod_{s=1}^{d} Q^s(\alpha, i_s) - \widetilde{f}_{i_1 \ldots i_d} \tag{7}$$

is the pointwise discrepancy, $f_{i_1 \ldots i_d} = \sum_{\alpha=1}^{r} \prod_{s=1}^{d} Q^s(\alpha, i_s)$ is the estimate of the tensor by canonical decomposition, $\widetilde{f}_{i_1 \ldots i_d}$ is the exact tensor (value of the function at nodes). Due to cumbersome and not lucid general expressions we shall provide in parallel some particular expressions corresponding considered six-dimensional case. In particular, we use the functional of discrepancy (6) in a form

$$\widetilde{\varepsilon}_{\Sigma} = 1/2(f_{ijklmp} - \widetilde{f}_{ijklmp}) \cdot (f_{ijklmp} - \widetilde{f}_{ijklmp}) / N^6, \tag{8}$$

where $f_{ijklmp} = \sum_{\alpha=1}^{r} Q^x(\alpha,i) \cdot Q^y(\alpha,j) \cdot Q^z(\alpha,k) \cdot Q^u(\alpha,l) \cdot Q^v(\alpha,m) \cdot Q^w(\alpha,p)$, $\widetilde{f}_{ijklmp}$ is an exact value of function at point $i,j,k,l,m,p$, $\widetilde{R}_{ijklmp} = f_{ijklmp} - \widetilde{f}_{ijklmp} = \sum_{\alpha=1}^{r} Q^x(\alpha,i) \cdot ... \cdot Q^w(\alpha,p) - \widetilde{f}_{ijklmp}$ is the pointwise discrepancy. The calculation of the functional (8) is laborious due to the high order of the tensor and we do not perform it. However, it is convenient to obtain for this functional some useful expressions that we shall use in more realistic approach. The gradient of discrepancy (6) has a form (we often use the notation $Q^s(\alpha,i_s) = Q^s_{\alpha,j_s}$):

$$\frac{\partial \widetilde{\varepsilon}_{\Sigma}}{\partial Q^c_{\alpha,j_c}} = \sum_{i_s, s=1...d, s\neq c} \widetilde{R}_{i_1...i_d} \cdot \prod_{s\neq c}^{d} Q^s(\alpha,i_s) / N^d. \tag{9}$$

In our case it is written as:

$$\nabla^x_{\alpha,i} \widetilde{\varepsilon}_{\Sigma} = \widetilde{R}_{ijklmp} \cdot \{Q^y(\alpha,j) \cdot ... \cdot Q^w(\alpha,p)\} / N^6,$$
$$...$$
$$\nabla^w_{\alpha,p} \widetilde{\varepsilon}_{\Sigma} = \widetilde{R}_{ijklmp} \cdot \{Q^x(\alpha,i) \cdot ... \cdot Q^v(\varepsilon,m)\} / N^6, \tag{10}$$

and is obtained from the variation of the discrepancy functional (8):

$$\Delta \widetilde{\varepsilon}_{\Sigma}(\Delta Q^x) = \{\sum_{\beta=1}^{r} Q^x(\beta,i) \cdot ... \cdot Q^w(\beta,p) - \widetilde{f}_{ijklmp}\} \cdot \{\sum_{\alpha=1}^{r} \Delta Q^x(\alpha,i) \cdot ... \cdot Q^w(\alpha,p)\} / N^6 = (\nabla^x_{\alpha,i} \varepsilon_{\Sigma}, \Delta Q^x(\alpha,i))$$
$$...$$
$$\Delta \widetilde{\varepsilon}_{\Sigma}(\Delta Q^w) = \{\sum_{\beta=1}^{r} Q^x(\beta,i) \cdot ... \cdot Q^w(\beta,p) - \widetilde{f}_{ijklmp}\} \cdot \{\sum_{\alpha=1}^{r} Q^x(\alpha,i) \cdot ... \cdot \Delta Q^w(\alpha,p)\} / N^6 = (\nabla^w_{\alpha,p} \varepsilon_{\Sigma}, \Delta Q^w(\alpha,p)). \tag{11}$$

The optimization search for the total set of the canonical decomposition cores

$$Q^c(\alpha,i) = \arg\min_{Q^c(\alpha,i)} \widetilde{\varepsilon}_{\Sigma}(Q^c(\alpha,i)), c = 1...d; i = 1...N_d; \alpha = 1...r \tag{12}$$

formally may be conducted using steepest descent iterations in all space of cores

$$\{Q^c(\alpha,i)\}^{n+1} = \{Q^c(\alpha,i)\}^n - \tau \nabla^c_{\alpha,i} \widetilde{\varepsilon}_{\Sigma} \tag{13}$$

along the gradient of the global discrepancy $\nabla^c_{\alpha,i} \widetilde{\varepsilon}_{\Sigma} \in R^{d \cdot r \cdot N}$ with the unique step $\tau$.

The Newton method may serve as another optimization option. It requires construction and inversion of the Hessian $H^{\Sigma}_{kc i_k i_c \gamma \beta} = \frac{\partial^2 \widetilde{\varepsilon}_{\Sigma}}{\partial Q^k(\gamma,i_k)\partial Q^c(\beta,i_c)} \in R^{d \times d \times N \times N \times r \times r}$ .

In practice, the application of the steepest descent (13) in space of the dimension $d \times N \times r$ is difficult due to ravine shape of the functional (2) caused by the mutual compensation of cores (relations of the form $\widetilde{Q}^1 \cdot \widetilde{Q}^2 = aQ^1 \cdot Q^2 / a$). The computational difficulties for the Newton method are caused by the need for the inversion of the Hessian of high dimension ($(d \times N \times r)^2$).

Below we consider the feasibility to reduce the dimension of the space, in which we perform an optimization, by $r$, and the dimension of the Hessian by $r^2$ that means $d \times N$ separate independent calculations in all nodes of cores.

## 3. RELATION OF GLOBAL DISCREPANCY EXTREMUM AND STATIONARY POINTS OF THE GRID OF LOCAL DISCREPANCIES

The global discrepancy (6) one may write as a sum of discrepancies on hyperplanes orthogonal to one of directions (here $c$), as

$$\widetilde{\varepsilon}_{\Sigma} = \sum_{i_c=1}^{N} \widetilde{\varepsilon}_{c,i_c} / N, \tag{14}$$

where local discrepancy contains the summation over the hyperplane

$$\widetilde{\varepsilon}_{c,i_c} = \sum_{\substack{i_s, s=1\ldots d \\ s \neq c}} \widetilde{R}_{i_1 \ldots i_d} \cdot \widetilde{R}_{i_1 \ldots i_d} / (2 \cdot N^{d-1}). \tag{15}$$

The grid of $d \times N$ local discrepancies $\widetilde{\varepsilon}_{c,i_c}$ is of special interest for us and is considered in details. The derivative of local discrepancy has a form

$$\frac{\partial \widetilde{\varepsilon}_{c,i_c}}{\partial Q^c_{\alpha,i_c}} = \sum_{i_s, s=1\ldots d, s \neq c} \widetilde{R}_{i_1 \ldots i_d} \cdot \prod_{s \neq c}^{d} Q^s(\alpha, i_s) / N^{d-1}. \tag{16}$$

In the case, for example, $x(i)$ ($\widetilde{\varepsilon}_{\Sigma} = \sum_{i=1}^{N} \widetilde{\varepsilon}_{x,i} / N$) the local discrepancy at point $i$ contains summation over hyperplane:

$$\widetilde{\varepsilon}_{x,i} = \sum_{jklmp} \widetilde{R}_{ijklmp} \cdot \widetilde{R}_{ijklmp} / (2 \cdot N^5). \tag{17}$$

Accordingly, the following relation holds between the components of the global discrepancy gradient and the derivatives of discrepancies on hyperplanes

$$\frac{\partial \widetilde{\varepsilon}_{\Sigma}}{\partial Q^x_{\alpha,i}} = \sum_{jklmp} \widetilde{R}_{ijklmp} \cdot [Q^y(\alpha, j) \cdot \ldots \cdot Q^w(\alpha, m)] / N^6 = \frac{\partial \widetilde{\varepsilon}_{x,i}}{\partial Q^x_{\alpha,i}} / N. \tag{18}$$

In general case of the discrepancy on hyperplane orthogonal to the direction $c$ the component of the gradient of the global functional $\widetilde{\varepsilon}_{\Sigma}$ on $Q^c_{i,\beta}$ is equal to the derivative of the local functional $\widetilde{\varepsilon}_{c,i}$ over this $Q^c_{i,\beta}$:

$$\frac{\partial \widetilde{\varepsilon}_{\Sigma}}{\partial Q^c_{\alpha,i_c}} = \sum_{i_s, s=1\ldots d, s \neq c} \widetilde{R}_{i_1 \ldots i_d} \cdot \prod_{s \neq c}^{d} Q^s(\alpha, i_s) / N^d = \frac{\partial \widetilde{\varepsilon}_{c,i_c}}{\partial Q^c_{\alpha,i_c}} / N. \tag{19}$$

From (19) one may see that the relation $\partial \widetilde{\varepsilon}_{.,i} / \partial Q_{\alpha,i}^c = 0$ entails $\partial \widetilde{\varepsilon}_\Sigma / \partial Q_{\alpha,i}^c = 0$, which corresponds the condition of the extremum of $\widetilde{\varepsilon}_\Sigma$. This causes

***Theorem 1***: ***The set of the canonical decomposition cores*** $Q_{\alpha,i_c}^c$ , ***for which the conditions для*** $\partial \widetilde{\varepsilon}_{c,i} / \partial Q_{\alpha,i}^c = 0$ ***are valid on the grid of*** $d \times N$ ***functionals*** $\widetilde{\varepsilon}_{c,i}$ , ***determines the point of extremum of global functional*** $\widetilde{\varepsilon}_\Sigma$ ***in space*** $Q_{\alpha,i_c}^c$ .

In result, instead the functional $\widetilde{\varepsilon}_\Sigma$ in the space of dimension $d \times N \times r$ we analyze the grid of $d \times N$ functionals $\widetilde{\varepsilon}_{c,i}$ each in the space of dimension $r$. As it takes place, we may seek for stationary points $\partial \widetilde{\varepsilon}_{c,i} / \partial Q_{\alpha,i}^c = 0$ on the grid by Newton Method

$$(Q^c(\alpha,i_c))^{n+1} = (Q^c(\alpha,i_c))^n - H_{(c,i_c)\alpha\gamma}^{-1} \cdot \partial \widetilde{\varepsilon}_{c,i_c} / \partial Q^c(\gamma,i_c) \qquad (20)$$

using the inversion of the local Hessian $H_{(c,i_c)\alpha\gamma} = \dfrac{\partial^2 \widetilde{\varepsilon}_{c,i_c}}{\partial Q^c(\alpha,i_c)\partial Q^c(\gamma,i_c)}$ . The set of stationary points for $d \times N$ functionals determines the point of extremum for $\widetilde{\varepsilon}_\Sigma(Q_{i,\beta}^c)$ .

The expression (20) is close in structure to the algorithm marked in [11] by ALS and, in our notations, expressed as

$$(Q^c(\beta,i_c))^{n+1} = H_{(c,i_c)\beta\gamma}^{-1} \cdot \varphi(\gamma,i_c) , \qquad (21)$$

where

$$\varphi(\gamma,i_c) = \sum_{i_s, s \neq c} \prod_{s=1, s \neq c}^{d} (Q^s(\gamma,i_s))^n \cdot \widetilde{f}_{i_1 \dots i_d} . \qquad (22)$$

In the strict sense, if we determined by calculations that the extremum of $\widetilde{\varepsilon}_\Sigma$ corresponds minimum, then all $\widetilde{\varepsilon}_{c,i}$ are minimums also, since if some $\widetilde{\varepsilon}_{c,i}$ is not minimum, then the sum $\widetilde{\varepsilon}_\Sigma = \sum_{i=1}^{N} \widetilde{\varepsilon}_{c,i} / N = \sum_{i,c}^{N,d} \widetilde{\varepsilon}_{c,i} /(d \cdot N)$ is not a minimum also. Thus we can seek points $\widetilde{\varepsilon}_{c,i}$, such that $\partial \widetilde{\varepsilon}_{c,i} / \partial Q_{\alpha,i}^c = 0$, using iterations of the gradient descent as in [18]

$$\{Q^c(\alpha,i_c)\}^{n+1} = \{Q^c(\alpha,i_c)\}^n - \tau_{c,i} \nabla_{c,i_c,\alpha} \widetilde{\varepsilon}_{c,i_c}, c = 1 \dots d; i_c = 1 \dots N . \qquad (23)$$

In general, the *Theorem 1* enables to significantly simplify the search for the canonical decomposition cores that provides the minimum of $\widetilde{\varepsilon}_\Sigma$. Instead the optimization step in the space $R^{r \times d \times N}$ , $d \times N$ optimizations in the space $R^r$ are performed.

Nevertheless, this approach remains difficult for tensors of the high order due to summation over tensor nodes on hyperplanes at calculation of derivatives. In order to further relax the computational limitations we shall apply the Monte-Carlo method at estimation of the integrals on hyperplanes.

## 4. CALCULATION OF LOCAL DISCREPANCIES AND THEIR DERIVATIVES USING MONTE-CARLO METHOD

The calculation of functionals (6,8), (15) and derivatives (17), (18) for the considered problems is very expensive from the computation time viewpoint due to high order of tensor. By this reason we estimated the approximation error norm (8) by the Monte-Carlo method

$$\varepsilon_{MC} = \sum_{e=1}^{L_{ens,t}} \{\sum_{\alpha=1}^{r} Q^x(\alpha,i_e) \cdot Q^y(\alpha,j_e) \cdot Q^z(\alpha,k_e) \cdot Q^t(\alpha,l_e) \cdot Q^y(\alpha,m_e) \cdot Q^v(\alpha,p_e) - \widetilde{f}_{i_e,j_e,k_e,l_e,m_e,p_e}\}^2 /(2L_{ens,t}), \quad (24)$$

where at every step of summation each index $i_e, j_e, k_e, l_e, m_e, p_e$ was chosen as a random and uniformly distributed from 1 to $N$ number. The amount of points of the ensemble $L_{ens,t} = 10^5 \div 10^6$ provided acceptable accuracy at moderate computation time (in $R^6$ at $N = 10^2$ the calculation of (8) requires summation over $10^{12}$ points).

The local discrepancy on hyperplane (15) similarly to (24) may be approximated using Monte-Carlo method as:

$$\varepsilon_{c,i_e} = \sum_{e=1,i_s^e=i_c,s=c}^{L_{ens}} R_{i_1^e...i_s^e...i_d^e} \cdot R_{i_1^e...i_s^e...i_d^e} /(2 \cdot L_{ens}), \quad (25)$$

where $L_{ens}$ is the number of the ensemble points on hyperplane. We use herein pointwise discrepancy

$$R_{i_1^e...i_s^e...i_d^e} = \sum_{\alpha=1}^{r} \prod_{s=1}^{d} Q^x(\alpha,i_s^e) - \widetilde{f}_{i_1^e...i_s^e...i_d^e}. \quad (26)$$

Let's consider for illustration one of local discrepancies (calculation of (17) using Monte-Carlo)

$$\varepsilon_{x,i} = 1/(2L_{ens}) \sum_{e=1}^{L_{ens}} \{(\sum_{\alpha} Q^x(\alpha,i) \cdot Q^y(\alpha,j_e^i) \cdot Q^z(\alpha,k_e^i)... - \widetilde{f}_{ij_e^i k_e^i l_e^i m_e^i p_e^i}) \cdot$$
$$\cdot (\sum_{\beta} Q^x(\beta,i) \cdot Q^y(\beta,j_e^i) \cdot Q^z(\beta,k_e^i)... - \widetilde{f}_{ij_e^i k_e^i l_e^i m_e^i p_e^i})\}, \quad (27)$$

where $\widetilde{f}_{ij_e^i k_e^i l_e^i m_e^i p_e^i}$ is the exact value of tensor at point $i, j_e^i, k_e^i, l_e^i, m_e^i, p_e^i$. The ensemble points $j_e^i, k_e^i, l_e^i, m_e^i, p_e^i$ are selected by indices $i$ and $e = 1...L_{ens}$. For $Q^x(\alpha,i)$ they are located on the hyperplane, orthogonal to $x$ at point $i$. So, at estimation of core $Q^x(\alpha,i)$ at point $i$, the discrepancy $\varepsilon_{x,i}$ depends on $L_{ens}$ points of the ensemble (Fig. 1), chosen randomly.

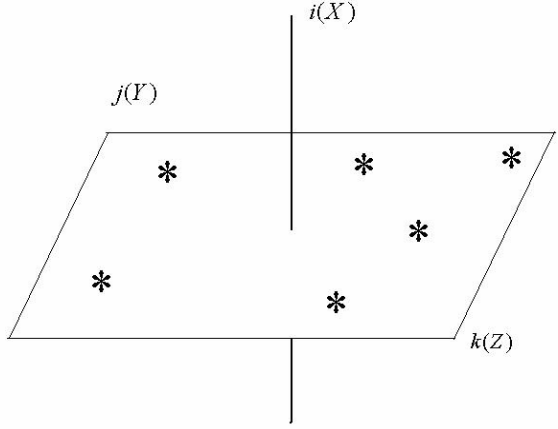

Fig. 1. Ensemble of randomly chosen points used for calculation of discrepancy $\varepsilon_{x,i}$

The grid of $d \times N$ local discrepancies (15) is of special interest for us since it enables to simplify significantly the algorithm. Similarly to (14), the global discrepancy may be considered as a sum of local discrepancies over hyperplanes:

$$\varepsilon_\Sigma = \sum_{i=1}^{N} \varepsilon_{x,i} / N = \sum_{c=1}^{d} \sum_{i=1}^{N} \varepsilon_{c,i} /(d \cdot N) \approx \widetilde{\varepsilon}_\Sigma \approx \varepsilon_{MC} . \tag{28}$$

The relations between derivatives of the global and local discrepancies, are valid similarly to (18)

$$\frac{\partial \varepsilon_\Sigma}{\partial Q_{\alpha,i}^x} = \sum_{e=1}^{L_{ens}} \{R_{ij_ek_el_em_ep_e} \cdot [Q^y(\alpha, j_e) \cdot ... \cdot Q^w(\alpha, p_e)]\} / L_{ens} = \frac{\partial \varepsilon_{x,i}}{\partial Q_{\alpha,i}} / N \tag{29}$$

and (19)

$$\frac{\partial \varepsilon_\Sigma}{\partial Q_{\alpha,i_c}^c} = \sum_{e=1,i_s^e=i_c,s=c}^{L_{ens}} R_{i_1^e...i_s^e...i_d^e} \cdot \prod_{s \neq c}^{d} Q^s(\alpha, i_s^e)\} / L_{ens} = \frac{\partial \varepsilon_{c,i_c}}{\partial Q_{\alpha,i_c}^c} / N . \tag{30}$$

Expressions (29), (30) differ from (18), (19) only by the calculation of integrals (convolutions) on hyperplanes using Monte-Carlo method. In other features, the further analysis is similar to expressions (6)-(23).

Let's consider the derivatives of functionals on hyperplanes in more details. For example, we disturb the core $Q^x(\beta, i)$ by $\Delta Q^x(\beta, i)$. The corresponding disturbance of $\varepsilon_{x,i}$ has a form (similar (11)):

$$\Delta \varepsilon_{x,i} = \sum_{e=1}^{L_{ens}} \{(\sum_{\beta=1}^{r} Q^x(\beta, i) \cdot Q^y(\beta, j_e^i) \cdot Q^z(\beta, k_e^i)... - \widetilde{f}_{ij_e^ik_e^il_e^im_e^ip_e^i}) \cdot$$
$$\cdot (\sum_{\alpha} \Delta Q^x(\alpha, i) \cdot Q^y(\alpha, j_e^i) \cdot Q^z(\alpha, k_e^i)... - \widetilde{f}_{ij_e^ik_e^il_e^im_e^ip_e^i})\} / L_{ens} . \tag{31}$$

The derivative of $\varepsilon_{x,i}$ over core $Q^x$ has a form:

$$\partial \varepsilon_{x,i} / \partial Q^x(\alpha,i) = \sum_{e=1}^{L_{ens}} \{(\sum_{\beta=1}^{r} Q^x(\beta,i) \cdot Q^y(\beta,j_e^i) \cdot Q^z(\beta,k_e^i)... - \widetilde{f}_{ij_e^i k_e^i l_e^i m_e^i p_e^i}) \cdot$$

$$\cdot Q^y(\alpha,j_e^i) \cdot Q^z(\alpha,k_e^i)...\} / L_{ens} . \tag{32}$$

Derivatives (29), (32) have $r$ component and requires the memory $\sim d \times N \times r$, the time cost of their calculation $\sim d \times N \times r \times L_{ens}$, they are local for every functional. All these features makes their application very attractive from the computational viewpoint.

## 5. NEWTON METHOD FOR THE SET OF DISCREPANCIES ON HYPERPLANES

The Newton method in many papers [13,15,14] is used for the search of the global discrepancy minimum. Usually, the expression for a Hessian is based on the Khatri-Rao product, that causes significant (and often nonrealistic) requirements to the computer memory and the time of computation.

Alternatives, not related with Khatri-Rao product, are presented in [11], where, local Hessians of dimension $(d \times r^2)^2$ are used at a minimization of linearized discrepancy by the Newton method. The method marked in [11] by ALS and implied the inversion for Hessians of dimension $r^2$ is the best from our opinion option.

In the frame of our approach (Theorem 1) we determine the extremum of the global discrepancy functional by the stationary points of the grid of local functionals. The Newton method is a natural algorithm for the search of stationary points. We use the gradients of local discrepancies and corresponding Hessians with dimension $r^2$. Let's consider this method at the example of the core $Q^x(\beta,i)$. The differentiation of gradient (32) provides expression for Hessian:

$$H_{(x,i)\beta\gamma} = \sum_{e=1}^{L_{ens}} \{[Q^y(\gamma,j_e^i) \cdot Q^z(\gamma,k_e^i)...] \cdot [Q^y(\beta,j_e^i) \cdot Q^z(\beta,k_e^i)...]\} / L_{ens} . \tag{33}$$

In general case:

$$H_{(c,i_c)\beta\gamma} = \frac{\partial}{\partial Q^c_{\beta,i_c}} \frac{\partial \varepsilon_{c,i_c}}{\partial Q^c_{\gamma,i_c}} = \sum_{e=1}^{L_{ens}} \prod_{s \neq c}^{d} Q^s(\beta,i_s^{e,i_c}) \cdot \prod_{s \neq c}^{d} Q^s(\gamma,i_s^{e,i_c}) / L_{ens} \tag{34}$$

Usually Hessian of discrepancy consist of two parts $H_{(x,i)\beta\gamma} = \frac{\partial^2 \varepsilon_{x,i}}{\partial Q^x(\beta,i)\partial Q^x(\gamma,i)} = \frac{1}{L_{ens}} \sum_{e=1}^{L_{ens}} \left( \frac{\partial f}{\partial Q^x(\beta,i)} \frac{\partial f}{\partial Q^x(\gamma,i)} + \frac{\partial^2 f}{\partial Q^x(\beta,i)\partial Q^x(\gamma,i)}(f - \widetilde{f}) \right)$. However, for the canonical decomposition (and the tensor train) second derivatives are equal to zero due to linear approximation of $f$ by every core. So, instead of complete Hessian we deals with Gauss-Newton matrix.

The degeneration of Hessian (34) at zero cores is the drawback from the computational viewpoint. It may be compensated by the Tikhonov regularization of zero order in the form $H_{(x,i)\beta\gamma} + \eta \delta_{\beta\gamma}$.

The availability of the easy computable Hessian enables effective estimation of the cores of canonical decomposition using Newton method

$$(Q^c(\alpha, i_c))^{n+1} = (Q^c(\alpha, i_c))^n - H^{-1}_{(c,i_c)\alpha\gamma} \cdot \partial\varepsilon_{c,i_c} / \partial Q^c(\gamma, i_c) \,. \tag{35}$$

The structure of corresponding algorithm at iteration $n+1$ has the form:

| *Algorithm 1. The search for canonical decomposition cores by Newton method* | | | |
|---|---|---|---|
| **for** $c = 1...d$ **do** | ! | | *cycle over coordinates* |
| **for** $i = 1...N_k$ **do** | ! | | *cycle over nodes along core* |
| **for** $\beta = 1...r$ **do** | ! | | *cycle over core layers* |
| $(Q^c(\alpha, i_c))^{n+1} = (Q^c(\alpha, i_c))^n - H^{-1}_{(c,i_c)\alpha\gamma} \cdot \partial\varepsilon_{c,i_c} / \partial Q^c(\gamma, i_c)$ | | ! | *summation over* $\gamma$ |
| **end for** | | | |

Herein we perform pointwise (on the grid $d \times N$ over coordinates and nodes along coordinate (core)) optimization in the space of core layers with gradient $\partial\varepsilon_{c,i_c} / \partial Q^c(\alpha, i_c) \in R^r$ and Hessian $H_{(k,i)\beta\gamma} \in R^{r \times r}$. The inversion of Hessian was performed by Gauss-Jordan method [20]. The time for single iteration is proportional to $d \times N \times r^3 \times L_{ens}$ (implying that the inversion of $r \times r$ matrix requires $r^3$ operations).

## 6. THE STEEPEST DESCENT METHOD FOR THE SET OF DISCREPANCIES ON HYPERPLANES

As it is mentioned above, one may apply the steepest descent method (23) with calculation of derivatives (30) by Monte-Carlo method as in [18] for the search of cores corresponding $\partial\varepsilon_{c,i_c} / \partial Q^c(\alpha, i_c) = 0$.

The structure of the algorithm is defined by cycles over coordinates, nodes and layers of cores:

| *Algorithm 2. The search for canonical decomposition cores by steepest descent* | | |
|---|---|---|
| **for** $c = 1...d$ **do** | ! | *cycle over coordinates* |
| **for** i$=1...N_c$ **do** | ! | *cycle over nodes along core* |
| **for** $\beta = 1...r$ **do** | ! | *cycle over core layers* |
| $(Q^c(\alpha, i_c))^{n+1} = (Q^c(\alpha, i_c))^n - \tau_{c,i_c} \cdot \partial\varepsilon_{c,i_c} / \partial Q^c(\alpha, i_c)$ | | |
| **end for** | | |

The time for single iteration is proportional to $r \times d \times N \times L_{ens}$.

## 7. RANDOM START

If the initial approximation of cores is selected to be unique (for example $Q^c(\alpha, i_c) = 1$), all values of gradients for different $\alpha$ coincide automatically and optimization fails to start. So, the expression $Q^c(\alpha, i_c) = 1 + N(\sigma)$ was used as the initial approximation. The normally distributed random value with the dispersion $\sigma$ was used. Naturally, there is no convergence at $\sigma = 0$. The instabilities were observed at great $\sigma \approx 10$. The optimal value of dispersion is about $\sigma = 0.1$, which was used in all numerical tests.

## 8. REGULARIZATION

Since the estimation of canonical decomposition cores is ill-posed problem [1,2], we applied the zero order Tikhonov regularization [19] (both for the steepest descent and for Newton method)

$$\varepsilon_{c,i_c} + \eta \sum_{\alpha=1}^{r} Q^c(\alpha, i_c)^2 / 2 .$$ (36)

Herein $\eta$ is the regularization coefficient. If the steepest descent is used without the regularization, the growth of the cores is observed with the difference between adjacent core nodes about 2-3 orders of magnitude. The Newton method without regularization does not operates. The iterations successfully converge for coefficients in the range $\eta = 10^{-4} \div 10^{-5}$.

## 9. RESULSTS OF NUMERICAL TESTS

The results of numerical tests on the approximation by canonical decomposition of six order tensors, obtained as the values of six-dimensional function in nodes of rectangular grid, are presented. The grid contained $N = 10^2$ nodes on every coordinate with the total number of nodes $10^{12}$. The construction of canonical decomposition was stopped at $\varepsilon_{\Sigma} \leq \varepsilon_2$, $\varepsilon_2 = 10^{-5} \div 10^{-7}$, the regularization coefficient was $\eta = 10^{-5}$, the rank of canonical decomposition was $r = 20$, the results provided for the number of points in ensemble on hyperplane $L_{ens} = 10^3$.

The global optimization is reduced to $6N$ local optimizations of discrepancies computed on hyperplanes in considered tests. The Newton method, the steepest descent and the method, which is marked in [11] as ALS with our modifications including applying (34) for calculation of Hessian and replacement of (22) by

$$\varphi(\gamma, i_c) = \sum_{\substack{e=1, i_s^e = i_c \\ s=c}}^{L_{ens}} \prod_{\substack{s=1, \\ s \neq c}}^{d} (Q^s(\gamma, i_s))^n \cdot \widetilde{f}_{i_1^e \dots i_s^e \dots i_d^e} .$$ (37)

Below we mark it as ALS.

As the first test, we consider the approximation of widespread [11] function in six dimensional space:

$$f = 1.0 / ((x/5)^2 + (y/5)^2 + (z/5)^2 + (u/5)^2 + (v/5)^2 + (w/5)^2)^{1/2} .$$ (38)

The Newton method required 3 iterations (1.2 minutes on Intel Core I5, 2,67 gHz) to converge to $\varepsilon_{\Sigma} = 10^{-6}$. ALS also required 3 iterations, practically with the same time. The steepest descent converged to this level for 50 iterations (44.3 minutes). In general, despite greater time for single iteration, the Newton method enables significant reduction of computation time if compare with the steepest descent. It also provides greater accuracy.

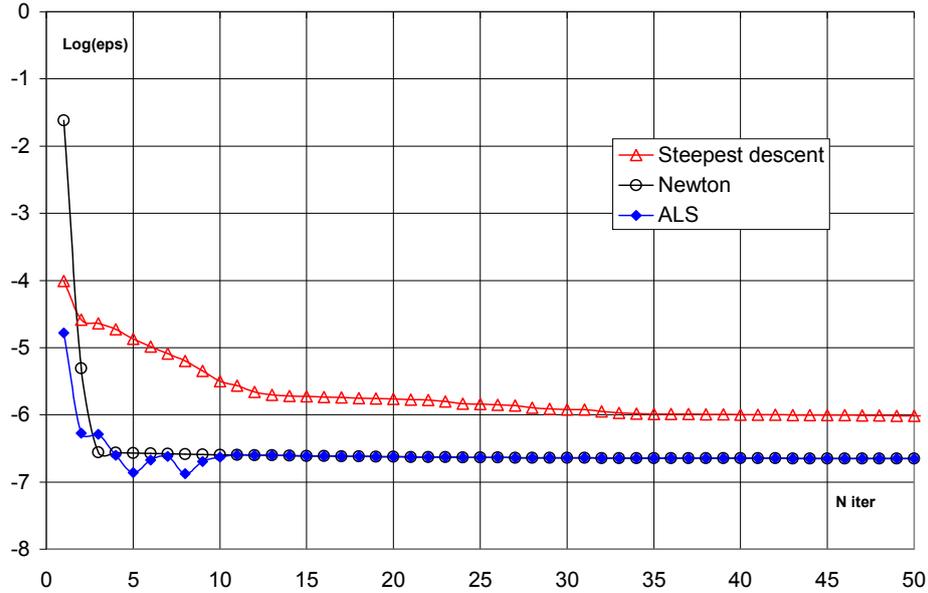

Fig. 2. The comparison of the convergence rate for steepest descent, Newton method and ALS for function (38).

Let's consider six-dimensional function defined by the sum of Gaussian and sines as a more difficult (due to complex non monotonic behavior) test

$$f = 5 \cdot \exp(-rad^2) + \sin(x/5) + \sin(y/5) + \sin(z/5) + \sin(u/5) + \sin(v/5) + \sin(w/5)$$
$$(rad = 0.001 \cdot ((ix - 50) + (iy - 50) + (iz - 50) + (iu - 50) + (iv - 50) + (iw - 50))) \,. \quad (39)$$

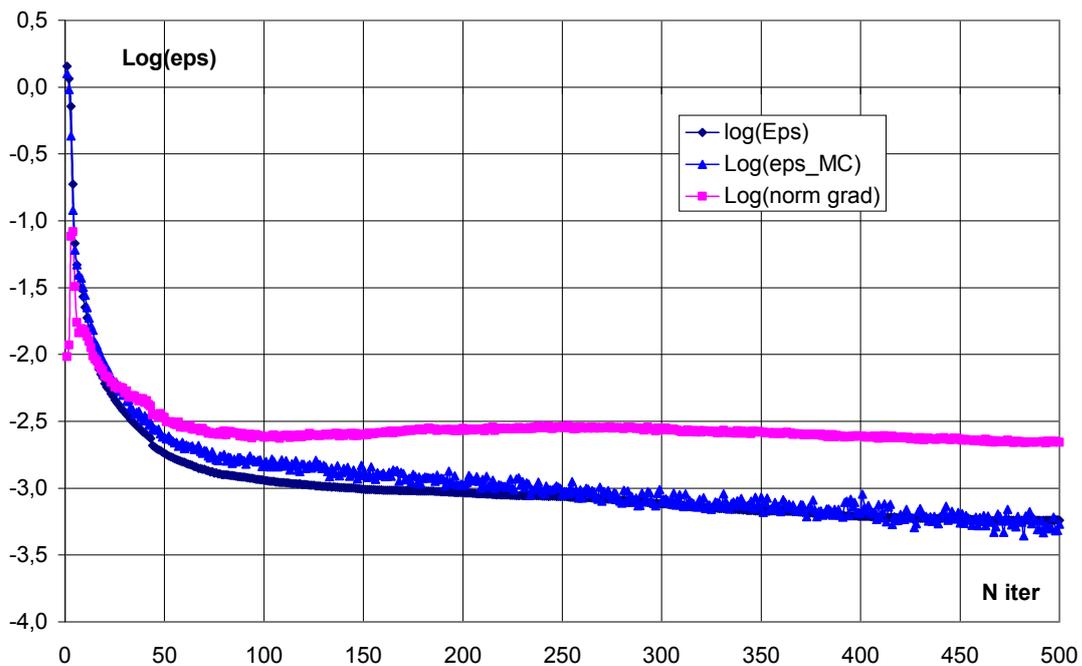

Fig. 3. Criteria of convergence in dependence on number of iterations of steepest descent for function (39)

Fig. 3 presents the behavior of different criteria of the convergence in dependence on the number of iterations on the example of steepest descent. The global discrepancy estimated via Monte-Carlo (24) $\varepsilon_{MC}$ at $L_{ens,t} = 10^5$ (Eps_MC), the global discrepancy $\varepsilon_\Sigma$, estimated by (26) (Eps) and the norm of the discrepancy gradient (grad norm) are presented in the logarithmic form. The discrepancies $\varepsilon_\Sigma$ (Eps) and $\varepsilon_{MC}$ (Eps_MC) demonstrate close magnitudes, the differences are related with the different structure and number of points in ensemble. The norm of the discrepancy gradient (grad norm) is not monotonic that hinder its application as the stopping criterion.

The comparison of the convergence rates for the steepest descent, Newton method and ALS is provided by Fig. 4. Global iterations stopped at $\varepsilon_\Sigma < 10^{-7}$. Newton method required 5.8 minutes. The steepest descent failed to achieve this accuracy. ALS also failed to get this accuracy, it required 26.4 minutes to achieve $\varepsilon_\Sigma < 10^{-6}$.

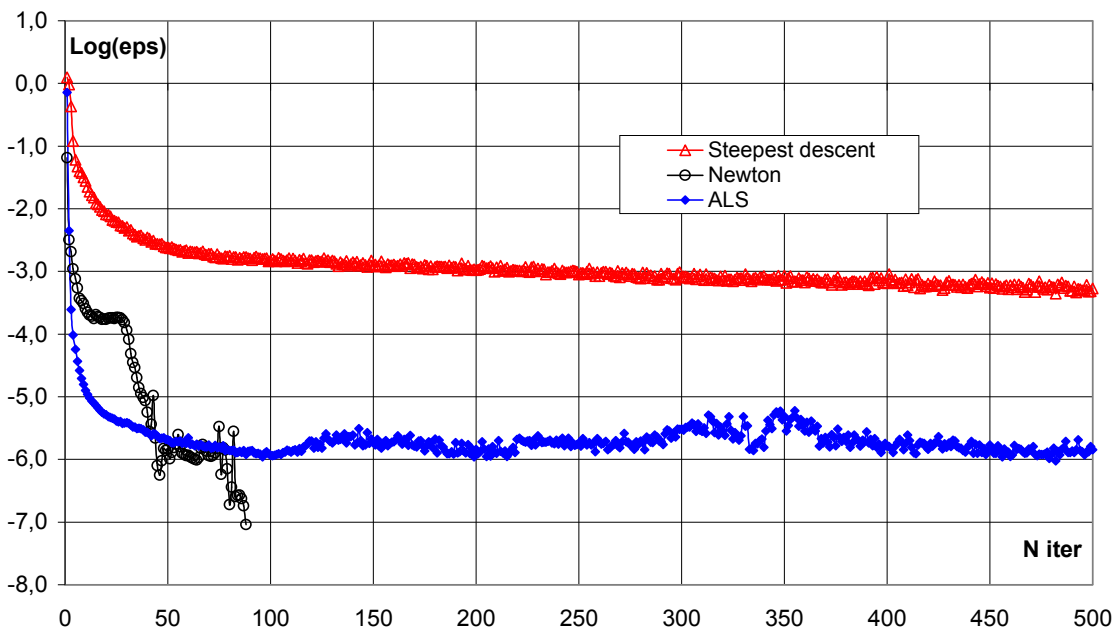

Fig. 4. The comparison of the convergence rates for the steepest descent, Newton method and ALS. Сравнение скорости сходимости наискорейшего спуска, метода Ньютона и ALS for function (39)

The section of function (39) by plane (x,y) centered over other coordinates and corresponding section of the approximation error are presented by Figs. 5 and 6.

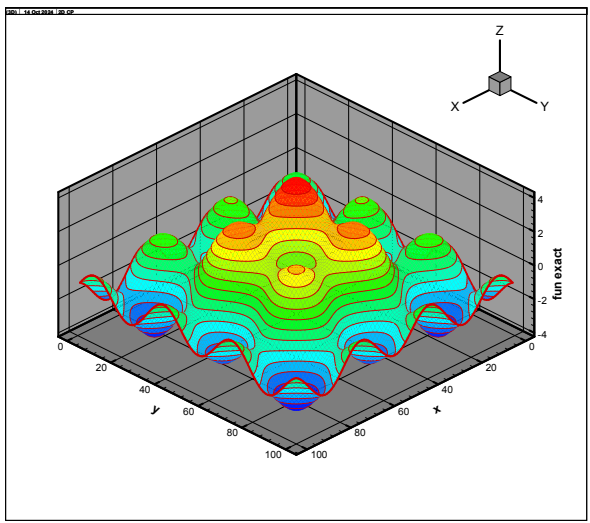

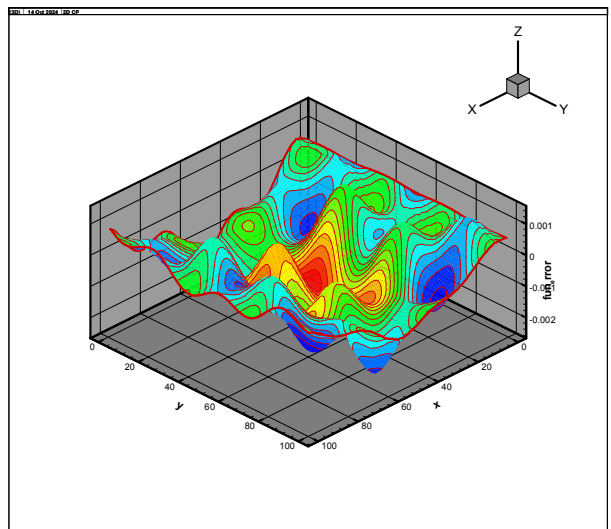

Fig. 5. Section of function (39) by plane (x,y)

Fig. 6. Section of approximation error for function (39) by plane (x,y)

Numerical tests demonstrated operability of the steepest descent, Newton method and ALS. In general, ALS is a bit worse Newton method from the accuracy and convergence rates viewpoints, however, the difference is not great and is visible not in all tests. The steepest descent method is inferior if compare with both Newton method and ALS.

The interest feature of the considered variant of the Newton method is related with the sparsity of results (some part of cores becomes zero) enabling the reduction of the decomposition rank. Neither steepest descent nor ALS demonstrates such behavior. The analysis of this peculiarity is far above the scope of present paper and will be presented in future papers.

## 10. DICSUSSION

Formally, the approximation of the tensor by the canonical decomposition eliminates the "curse of dimensionality" since the tensor has $N^d$ nodes while the canonical decomposition is described by $d \cdot N \cdot r$ parameters. Nevertheless, the approximation of tensor by the canonical decomposition is computationally expensive. The reason is that the majority of algorithms (known to authors) in some way uses the tensor (or its matricization without change of the number of nodes) and the Khatri-Rao product. The matrix, composed by the Khatri-Rao products contains $N^{d-1} + r$ elements [13]. The time of its inversion is proportional to $(N^{d-1} + r)^{3/2}$, and a time of single iteration is proportional $d \cdot (N^{d-1} + r)^{3/2}$. So, the application of the Khatri-Rao product is plagued by the "curse of dimensionality" both from the viewpoint memory and the time of computation and hardly can be used for tensors above third order.

In the considered variant of Newton method the local Hessian has $r^2$ параметров, correspondingly the time of its inversion $\sim r^3$, that is close to the algorithm marked in [11] by ALS. The total cycle of estimation of cores contains $d \times N$ local calculations at one global iteration. The time for computation of discrepancies (integrals) is proportionate $L_{ens}$. In result, the time of one global iteration is proportional to $d \times N \times r^3 \times L_{ens}$ and increase over $N$ more slow if compare with ALS by version of [11], where corresponding time is proportional to $d \times N \times r^3 \times N^{d-1}$.

In general, the necessity for a relative great ensemble forms the "neck of the bottle" for the considered approach from the viewpoint of computation time. If the ensemble of points on the hyperplane is chosen randomly (as in present tests) the error of calculation of integral by Monte-Carlo method does not depend on the dimension of space, but has rather low convergence rate over the ensemble size $\sim (L_{ens})^{-1/2}$. The acceleration of convergence rate is feasible by refusing from use a random uniform ensemble. Quasi Monte-Carlo [22], Latin Hyper cube [20, 23] and Sobol's sequences [21] have the error asymptotics close to $\sim (L_{ens})^{-1}$ and provide certain hope for acceleration of the approach.

The information in considered algorithms is transferred only via state of cores and there is no prohibition on the growth of discrepancy in total (over all cores) iteration, so certain non monotonicity (and even increase) of global discrepancy occurs, which, nevertheless, does not hinder the efficiency of algorithms.

The considered approach can be applied for the tensor train also, however, the detailed analysis of such possibility is far from the topic of current paper.

## CONCLUSIONS

The method for calculation of the canonical decomposition for high order tensor approximation is offered, which applies the Monte-Carlo method for calculation of multidimensional integrals and a minimization of the global norm of approximation error via the search for stationary points for the grid of discrepancies computed on hyperplanes.

The considered method does not use matricization of approximated tensor and direct calculation of multidimensional integrals that enables significant economy both the computer memory and the time of computation. It may be realized by the steepest descent, ALS, and Newton method, which implies an inversion of matrices of the dimension that is equal to the rank of decomposition.

The numerical tests for the tensors of six order, containing $10^{12}$ nodes, demonstrated the efficiency of this method and its ability to operate using office computers.

REFERENCES

1. *Hackbusch W.* Tensor Spaces and Numerical Tensor Calculus. Springer, 2012.
2. *Yorick H., Willi-Hans S.* Matrix Calculus, Kronecker Product and Tensor Product: A Practical Approach to Linear Algebra, Multilinear Algebra and Tensor Calculus with Software Implementation, Singapore: WSP Ltd., 2019.
3. *Oseledets I. V.* Tensor- train decomposition//SIAM J. Sci. Comput. 2011. V. 33. P. 2295–2317.
4. *G. Beylkin, M. Mohlenkamp.* Algorithms for Numerical Analysis in High Dimensions// SIAM Journal on Scientific Computing. 2005. V. 26 (6). P. 2133-2159.
5. *Kazeev V. A., Tyrtyshnikov E. E.* The structure of the Hessian and the efficient implementation of Newton's method in the problem of the canonical approximation of tensors//Comput. Math. Math. Phys., 2010. V**. 50**. N. 6. 927–945
6. *Oseledets I., Tyrtyshnikov E.* TT-cross approximation for multidimensional arrays// Linear Algebra Appl. 2010. V. 432. P. 70–88.
7. *Khatri C.G. and Rao C.R.* Solutions to some functional equations and their applications to characterization of probability distributions//Sankhya. 1968. Ser. A. *V.* 30. P. 167-180.
8. *Siaminou I., Liavas A. P.* An Accelerated Stochastic Gradient for Canonical Polyadic Decomposition// arXiv:2109.13964v1. 2021.
9. *Rabanser S., Shchur O., Gunnemann S.* Introduction to Tensor Decompositions and their Applications in Machine Learning//arXiv:1711.10781v1. 2017.
10. Uschmajew A. Local convergence of the alternating least squares algorithm for canonical tensor approximation//SIAM J. Matrix Anal. Appl. 2012. V. 33 (2). P. 639–652.
11. *Oseledets I. V., Savostyanov D. V.* Minimization methods for approximating tensors and their comparison//Comput. Math. Math. Phys. 2006. V. 46. N. 10. P. 1641–1650
12. *Acar E., Dunlavy D. M., and Kolda T. G.* A scalable optimization approach for fitting canonical tensor decompositions// J. Chemometrics. 2011. V. 25. P. 67-86.

13. Yanzhao Cao, Somak Das, Luke Oeding, and Hans-Werner Van Wyk, Analysis of the Stochastic Alternating Least Squares Method for the Decomposition of Random Tensors// arXiv:2004.12530v1. 2020.

14. *Singh N., Ma L., Yang H., and Solomonik E.* Comparison of Accuracy and Scalability of Gauss-Newton and Alternating Least Squares for CP Decompositions//arXiv:1910.12331v2. 2019.

15. *Phan A.-H., Tichavsky P., and Cichocki A.* Low complexity damped Gauss-Newton algorithms for CANDECOMP/PARAFAC//SIAM J. on Matrix Analysis and Applications. 2013. V. 34(1). P. 126-147.

16. *Nion D. and De Lathauwer L.* An enhanced line search scheme for complex-valued tensor decompositions. Application in DS-CDMA//Signal Processing. 2008. V. 88(3). P. 749-755.

17. *Rajih M., Comon P., and Harshman R.* A. Enhanced line search: A novel method to accelerate PARAFAC//SIAM J. Matrix Analysis and Appl. 2008. V. 30(3). P. 1128-1147.

18. *Alekseev A.K., Bondarev A.E., Pyatakova Yu. S.* On Application of Canonical Decomposition for the Visualization of Results of Multiparameter Computations//Scientific Visualization. 2023. V. 15, N. 4. P. 12 – 23

19. *Tikhonov A.N., Arsenin V.Y.* Solutions of Ill-Posed Problems. Winston and Sons, Washington DC 1977

20. *Press W. H., Flannery B. P., Teukolsky S. A., Vetterling W. T.* Numerical Recipes in Fortran 77: The Art of scientific computing, Cambridge Univ. Press, 1992

21. *Sobol I. M.* On the distribution of points in a cube and the approximate evaluation of integrals//U.S.S.R. Comput. Math. Math. Phys. 1967. V. **7.** N**.** 4.P. 86–112

22. *Hastings W. K.* Monte Carlo sampling methods using Markov chains and their applications//Biometrika. 1970. V. 57. P. 97–109.

23. *Viana F. A. C.* A Tutorial on Latin Hypercube Design of Experiments//Qual. Reliab. Engng. Int. 2016, V. 32. P. 1975– 1985.